\documentclass[12pt]{amsart}
\usepackage{amssymb}
\usepackage{amsmath}
\usepackage{pictex}
\usepackage{url}
\allowdisplaybreaks{}

\newcommand{\abs}[1]{\left\lvert #1 \right\rvert} 
\newcommand\Acal{\mathcal{A}}
\newcommand\AND{\quad\text{and}\quad}
\newcommand\CC{\mathfrak L}
\newcommand\Ccal{\mathcal C}
\newcommand\ed{\text{\rm bd}}
\newcommand\ef{[x,y]}
\newcommand\Ef{E}
\newcommand\Ex{\mathsf{E}}

\newcommand\GG{\mathfrak{G}}

\newcommand{\GA}{\GG{}\!_A}
\newcommand{\hatGA}{\widehat{\GG}\!_A}
\DeclareMathOperator{\GL}{\mathrm{GL}}

\newcommand\HH{{\mathfrak H}}

\newcommand\IC{\mathbb C}
\newcommand\id{\text{\it id}}
\newcommand\idG{e}
\newcommand\idH{o}
\newcommand\la{\lambda}
\newcommand\msf{\mathfrak{m}}
\newcommand\nnu{\overline{\nu}}
\newcommand{\norm}[1]{\left\lVert #1 \right\rVert} 
\newcommand\nsf{\mathfrak{n}}
\newcommand\psf{\mathsf{p}}
\newcommand\Prob{\mathsf{Pr}}
\newcommand\R{\mathbb R}
\newcommand\refl[1]{\check{#1}}
\newcommand\Rsf{R}
\newcommand\spec{\operatorname{\sf spec}}
\newcommand\supp{\operatorname{\sf supp}}
\newcommand\uno{\iota}
\newcommand\wt{\widetilde}
\newcommand\Z{\mathbb Z}

\DeclareMathOperator{\diag}{\mathsf{diag}}
\DeclareMathOperator{\linspan}{\mathsf{span}}
\DeclareMathOperator{\Tr}{\mathsf{Tr}}

\numberwithin{equation}{section}

\makeatletter
  \DeclareRobustCommand\em
         {\@nomath\em \ifdim \fontdimen\@ne\font >\z@
                       \upshape \else \slshape \fi}
\makeatother

\newtheoremstyle{mythm}
  {9pt}
  {9pt}
  {\slshape}
  {0pt}
  {\bfseries}
  {}
  { }
  {\thmnumber{(#2)}\thmname{ #1}\thmnote{ #3}}

\newtheoremstyle{mydef}
  {9pt}
  {9pt}
  {\normalfont}
  {0pt}
  {\bfseries}
  {}
  { }
  {\thmnumber{(#2)}\thmname{ #1}\thmnote{ #3}}

\theoremstyle{mythm}
\newtheorem{thm}[equation]{Theorem.}
\newtheorem{pro}[equation]{Proposition.}
\newtheorem{lem}[equation]{Lemma.}
\newtheorem{cor}[equation]{Corollary.}

\theoremstyle{mydef}
\newtheorem{dfn}[equation]{Definition.}

\begin{document}$\,$ \vspace{-1truecm}
\title[On the spectrum of lamplighter groups and percolation clusters]{On the
    spectrum of lamplighter groups and percolation clusters\\
    {\mdseries\small\footnotesize(Extended Version)}}
\author{Franz LEHNER, Markus NEUHAUSER and Wolfgang WOESS}
\address{\parbox{.8\linewidth}{Institut f\"ur Mathematische Strukturtheorie,
TU Graz,\\
Steyrergasse 30, 8010 Graz, Austria (F.L. \& W.W.) \\}}
\email{lehner@finanz.math.tu-graz.ac.at, woess@tugraz.at}
\address{\parbox{.8\linewidth}{Lehrstuhl A f\"ur Mathematik, RWTH Aachen,\\
52056 Aachen, Germany (M.N.)\\}}
\email{markus.neuhauser@matha.rwth-aachen.de}
\date{\today} 
\thanks{Partially supported by Austrian Science Fund (FWF) P18703-N18}
\begin{abstract}
Let $\GG$ be a finitely generated group and $X$ its Cayley graph with
respect to a finite, symmetric generating set $S$. Furthermore,
let $\HH$ be a finite group and $\HH \wr \GG$ the lamplighter 
group (wreath product) over $\GG$ with group of ``lamps'' $\HH$.
We show that the spectral measure (Plancherel measure) of any symmetric
``switch--walk--switch'' random walk on $\HH \wr \GG$ coincides
with the expected spectral measure (integrated density of states)
of the random walk with absorbing boundary on the cluster of the group 
identity for Bernoulli site percolation on $X$ with parameter 
$\psf = 1/|\HH|$. The return probabilities of the
lamplighter random walk coincide with the expected (annealed)
return probabilities on the percolation cluster. In particular, 
if the clusters of percolation with parameter $\psf$ are almost surely 
finite then the spectrum 
of the lamplighter group is pure point. This generalizes 
results of Grigorchuk and \.Zuk, resp.\ Dicks and Schick regarding
the case when $\GG$ is infinite cyclic.

Analogous results relate bond percolation with another lamplighter random
walk. In general, the integrated density of states of 
site (or bond) percolation with arbitrary parameter $\psf$ is always related
with the Plancherel measure of a convolution operator by a signed
measure on $\HH \wr \GG\,$, where $\HH = \Z$ or another suitable group.
\end{abstract}

\maketitle

\baselineskip 15pt


\section{Introduction}\label{sec:intro}

\subsection{Lamplighter random walks}

Let $\GG$ be a finitely generated group
and $\HH$ a finite group with unit elements $\idG$ and
$\idH$, respectively.
The \emph{wreath product} or \emph{lamplighter group} $\HH \wr \GG$ is the 
semidirect product $\CC \rtimes \GG$, where $\CC = \bigoplus_{\GG} \HH$
is the group of \emph{configurations} 
$\eta: \GG \to \HH$ with finite support 
$\supp(\eta) = \{ x \in \GG : \eta(x) \ne \idH \}\,$. The group operation in
$\CC$ is pointwise
multiplication in $\HH\,$, and its unit element $\uno$ is given by $\uno(x) = \idH$
for all $x \in \GG$. The left action of $\GG$ on $\CC$ is 
$
L_g\eta(x) = \eta(g^{-1}x)\,,
$ 
so that the group operation in $\HH \wr \GG$ is
$$
(\eta,g)(\eta',g') = (\eta \cdot L_g\eta',gg')\,.
$$ 
We can embed $\GG$ and $\HH$ into $\HH \wr \GG$ via the mappings
$$
g \mapsto (\uno, g) \AND h \mapsto (\eta_{\idG}^h,\idG)\,,\quad\text{where
for}\;g \in \GG\,,\;h \in\HH\,, 
\quad \eta_g^{h}(x) = \begin{cases} 
      h\,,&\text{if}\; x=g\,,\\
      \idH\,,&\text{otherwise.}
		      \end{cases}
$$			       
Now let $\mu$ be a symmetric probability measure on $\GG$ whose support
$S$ is finite and generates $\GG$. The \emph{random walk on $\GG$ with law}
$\mu$ is the Markov chain with transition probabilities 
$p(x,y)= \mu(x^{-1}y)$, $x, y \in \GG$.  
The \emph{Cayley graph} $X(\GG,S)$ of $\GG$ with respect to $S$ has vertex
set $\GG$ and the unoriented edges $[x,xs]$, where $x \in \GG$ and 
$s \in S$. The steps of our random walk follow the edges of this graph,
and the most natural case is the one where $\mu$ is equidistributed on
$S$, in which case it generates the \emph{simple random walk} on $X$.

Also, we let $\nu$ be equidistribution on $\HH$. Via the above embedding, we
consider $\mu$ and $\nu$ as probability measures on $\HH \wr \GG$, and build 
the convolution
\begin{equation}\label{eq:mu}
\wt \mu = \nu*\mu*\nu\,.
\end{equation}
This is a symmetric probability measure whose support 
$$
\supp(\wt\mu) = \{ (\eta_{\idG}^h\cdot\eta_s^{h'},s) : h, h' \in \HH,\;
s \in S\, \}
$$
generates the lamplighter group. (If $\mu$ is equidistributed on $S$ and
$\idG \notin S$ then $\wt\mu$ is also equidistributed on its support.) 
It gives rise to the \emph{switch-walk-switch}
lamplighter walk: there are lamps at the vertices of $X(\GG,S)$ whose possible 
states are encoded by the group $\HH$, and $\idH$ is the state ``off'' of a 
lamp. Initially, all lamps are off. A lamplighter performs simple random walk
on $X(\GG,S)$,  starting at $\idG$. At each step, (s)he first puts the lamp at the current position
to a random state, then makes a move in $X(\GG,S)$, and finally puts the lamp at
the new position to a random state. Each current configuration of the lamps 
plus the current position of the lamplighter is a pair 
$(\eta,g) \in \HH \wr \GG\,$. This process is the random walk with
law $\wt \mu$ on $\HH \wr \GG\,$.  

In this note, one of the objects that we are interested in is the spectrum of 
the transition operator of the lamplighter random walk, that is, the
right convolution operator $\Rsf_{\wt\mu}\,$, acting on functions 
$F \in \ell^2(\HH \wr \GG)$ by $F \mapsto F * \wt\mu$. It is self-adjoint. 
(If, more generally, $\wt\mu$ was not symmetric, we would have to convolve
by the reflection of $\wt\mu$.) 
More precisely, we are interested 
in the \emph{Plancherel measure},
also called \emph{Kesten--von Neumann--Serre} spectral measure
by some authors. This is the on--diagonal element of the resolution of
the identity of our convolution operator, or equivalently, the probability
measure $\msf$ on $\spec(\Rsf_{\mu}) \subset [-1\,,\,1]$ whose moments  
\begin{equation}\label{eq:lampmoments}
\int_{[-1\,,\,1]} t^n\,d\msf(t) = \wt\mu^{(n)}(\uno,\idG)\,,\quad n \ge 0\,,
\end{equation}
are the probabilities that the lamplighter is back to the starting point 
at step $n$ with
all lamps switched off. Here, $\wt\mu^{(n)}$ denotes the $n$-th convolution power
of $\wt\mu$. 

We remark that due to the choice of $\nu$, that spectrum depends only
on the cardinality and not on the specific group structure of $\HH\,$.

In the case when the base group $\GG$ is infinite cyclic,
{\sc Grigorchuk and \.Zuk}~\cite{GrZu} (for $|\HH|=2$) and 
{\sc Dicks and Schick}~\cite{DiSc} (for arbitrary finite $\HH$) have computed
this spectrum and the Plancherel measure explicitly for simple random walk; 
a more elementary explicit computation that applies to a larger class of 
graphs was done by {\sc Bartholdi and Woess}~\cite{BaWo}. To be precise, those 
computations apply to a slightly different variant of the lamplighter walk, 
but carry over immediately to the ``switch--walk--switch'' model. 
In those results, it turns out that the spectrum is \emph{pure point;} 
a complete orthonormal system of finitely supported eigenfunctions of
the operator can be computed explicitly. 

Pure point spectra are quite familiar in the context of fractal structures,
see e.g.\ {\sc Teplyaev}~\cite{Te}, {\sc Kr\"on}~\cite{Kr} and {\sc
Sabot}~\cite{Sa} for rigorous results.
For random walks on groups, the above was the first example of such a
phenomenon. One reason can be found in the fact that $\HH \wr \Z$ has an
inherent structure of self-similarity, compare with {\sc Bartholdi,
Grigorchuk and Nekrashevych}~\cite{BaGrNe}. Later, further classes of groups
and convolution operators on them with pure point spectrum were found
by {\sc Bartholdi, Neuhauser and Woess}~\cite{BaNeWo}. 
The origin of this paper was the question at which level of generality a
pure point spectrum occurs for typical lamplighter random walks on general
wreath products $\HH \wr \GG$.

\subsection{Percolation clusters}\label{ssec:percolationclusters}
Consider our Cayley graph $X(\GG,S)$ and let $0 < \psf < 1$. In \emph{Bernoulli
site percolation} with parameter $\psf$, we have i.i.d. Bernoulli random
variables $Y_x\,$, $x \in \GG\,$, sitting at the vertices of $X$, with
$$
\Prob_{\psf}[Y_x=1] = \psf \AND \Prob_{\psf}[Y_x=0] = 1-\psf\,.
$$
(The index $\psf$ refers to the parameter.) 
We can realize those random variables on a suitable probability space 
$(\Omega, \Acal, \Prob_{\psf})$. In the percolation process, for every 
$\omega \in \Omega$, we keep all vertices where $Y_x(\omega)=1$
(the \emph{open} ones) and delete (close) the other ones. 
We consider the resulting induced subgraph spanned by the open vertices:
its edges are those with both endpoints open, while the other edges
are deleted. This graph 
falls apart into connected components, which are random subgraphs of $X$. 
If $x$ is open, then its component $C(x)=C_{\omega}(x)$ consists of all open 
vertices that are connected to $x$ by a path in the graph $X$ all whose 
vertices are open. 
If $x$ is closed then we set $C_{\omega}(x) = \emptyset$. 
It is well known that
there is a \emph{critical parameter} $\psf_c$ such that for any vertex~$x$,
$$
\Prob_{\psf}[C(x)\;\text{is finite}] = 1\,,\;\;\text{if}\;\;\psf < \psf_c\,,\AND
\Prob_{\psf}[C(x)\;\text{is finite}] < 1\,,\;\;\text{if}\;\;\psf > \psf_c\,.
$$
This property is independent of the specific vertex $x$. 
One usually writes 
$$
\theta(\psf) = \Prob_{\psf}[C(x)\;\text{is infinite}]\,,
$$
so that $\theta(\psf) = 0$ for $\psf < \psf_c$ and $\theta(\psf) > 0$
for $\psf < \psf_c\,$.
The behaviour at $\psf=\psf_c$ is more delicate. One always has
$\psf_c > 0$. The standard monograph about percolation on the integer
lattices is the one of {\sc Grimmet}~\cite{Gr}. A systematic study of
percolation on general Cayley graphs, as well as vertex-transitive graphs,
was initiated in the 1990s, see e.g. {\sc Benjamini, Lyons, Peres and
Schramm}~\cite{BLPS}. A quite complete account is given in the forthcoming
book by {\sc Lyons with Peres}~\cite{LyPe}.

Now we consider the restriction of the simple random walk on $X$ to 
the random graph $C_{\omega}= C_{\omega}(\idG)\,$, when $\idG$ is open. 

In general, for any finite or infinite, connected subgraph $A$
of $X$ that contains $\idG$, we define the transition probabilities
\begin{equation}\label{eq:clusterRW}
p_A(x,y) = \begin{cases} p(x,y)=\mu(x^{-1}y)\,,&\text{if both}\; 
                            x, y \in A\,,\\
                       0\,,&\text{otherwise.}
         \end{cases}
\end{equation}
The transition matrix $P_A$ is strictly substochastic at the points 
of $A$ which have a neighbour in the complement.
That is, there is a positive probability that the random walk is absorbed 
(dies) at such a point. The $n$-step transition probability
$p_A^{(n)}(x,y)$ can be interpreted as the probability that the 
simple random walk on $X$ moves from $x$ to $y$ in $n$ steps before
leaving $A$. With this interpretation, it also makes
sense to set $p_A^{(0)}(\idG,\idG)=1$ in the degenerate case 
when $A = \emptyset\,$. 
Now $P_A$ acts on functions $f \in \ell^2(C)$, as well as 
$f \in \ell^2(\GG)$, in the usual way by 
$P_Af(x) = \sum_y p_A(x,y)f(y)$. (Contrary to $P=R_{\mu}\,$, this is of course
not a right convolution operator.)
We also admit $A = \emptyset$, with
$P_{\emptyset}f = 0$.  For any $A$, the operator $P_A$ is self-adjoint, and we
can again consider the  diagonal element at $\idG$ of its
spectral resolution. This is the probability measure $\nsf_A$
on $\spec(P_A) \subset [-1\,,\,1]$ whose moments are
\begin{equation}\label{eq:clustermoments}
\int_{[-1\,,\,1]} t^n\,d\nsf_A(t) = p_A^{(n)}(\idG,\idG)\,,\quad n \ge 0\,.
\end{equation}
We are interested in random transition operator $P_{\omega} = P_{C_{\omega}}\,$,
the associated random spectral measure $\nsf_{\omega} = \nsf_{C_{\omega}}$
and its moments $p_{\omega}^{(n)}(\idG,\idG) = p_{C_{\omega}}^{(n)}(\idG,\idG)$.

\subsection{Main results}

There is a surprisingly simple relation between the lamplighter random walk
and the random walk with absorbing boundary on the percolation cluster.
The following statements refer to the setup described above. 

\begin{thm}\label{thm:spec-meas}
Consider site percolation on $X(\GG,S)$ with parameter $\psf=1/|\HH|$. 
The Plancherel measure $\msf$ of the lamplighter random walk on 
$\HH \wr \GG$ with law $\wt\mu$ given by \eqref{eq:mu} and the spectral measures
$\nsf_{\omega}$ of the random walk with
absorbing boundary given by \eqref{eq:clusterRW} on the cluster of $\idG$
are related by
\begin{equation}\label{eq:meas-ident}
\msf(B) = \Ex_\psf \bigl( \nsf_{\omega}(B) \bigr) \quad 
\text{for every Borel set}\;\;
B \subset \R\,.
\end{equation}
Equivalently,
the return probabilities of the respective random walks satisfy
\begin{equation}\label{eq:prob-ident}
\wt\mu^{(n)}(\uno,\idG) = \Ex_{\psf} \bigl( p_{\omega}^{(n)}(\idG,\idG) \bigr)
\quad \text{for all}\; n \ge 0\,.
\end{equation}
\end{thm}

Here, $\Ex_{\psf}$ refers of course to expectation on 
$(\Omega, \Acal, \Prob_{\psf})$.

In the mathematical physics literature, the measure on the right hand
side of \eqref{eq:meas-ident} is sometimes called the 
\emph{integrated density of states,} see e.g.\ 
{\sc Kirsch and M\"uller}~\cite{KiMu} and {\sc Veseli\'c}~\cite{Ve}. 
In the terminology of random walk in random environment, the
expected return probabilities on the right hand side of \eqref{eq:prob-ident}
are often called \emph{annealed} return probabilities. 

We remark that the simple identity \eqref{eq:prob-ident} holds for an arbitrary 
probability measure $\mu$ on $\GG$ in the place of one that is finitely 
supported and symmetric. 

\begin{dfn}\label{def:C} Let $\Ccal$ be the collection of all finite connected subgraphs of 
$X(\GG,S)$ that contain $\idG\,$, plus the empty set. 
\end{dfn}

Regarding the spectrum,  we have the following.
Recall that the \emph{point spectrum} of our operator~$R_{\wt\mu}$ is the set of eigenvalues, i.e., the set 
$$
\{\lambda\in\IC : \text{ $\lambda I-R_{\wt\mu}$ is not injective} \}.
$$

\begin{thm}\label{thm:pointspec} 
The point spectrum $\spec(\Rsf_{\wt\mu})$ of the right convolution operator on 
the space $\ell^2(\HH \wr \GG)$ by the measure $\wt\mu$ of \eqref{eq:mu} comprises the set
$$
\Lambda = 
\bigcup \bigl\{ \spec(P_A) : A \in \Ccal\,,\;A \ne \emptyset \bigr\}
\cup \{ 0 \}\,.
$$
For each eigenvalue $\lambda \in \Lambda$, the eigenspace contains infinitely
many linearly independent eigenfunctions with finite support.
\end{thm}

\begin{thm}\label{thm:pure} If $\psf=1/|\HH|$ is such that site percolation
on $X(\GG,S)$ satisfies $\theta(\psf) = 0$, that is,
$$
\Prob_{\psf}[C(\idG)\;\text{is finite}] = \sum_{A\in\Ccal}\Prob_{\psf}[C(\idG)=A]= 1\,,
$$
then $\spec(\Rsf_{\wt\mu})$ is pure point; it is the closure of $\Lambda$. 
There is a complete orthonormal system in $\ell^2(\HH \wr \GG)$ consisting
of finitely supported eigenfunctions associated with the eigenvalues in 
$\Lambda$.
\end{thm}

The special case $\GG=\Z$ was treated
in \cite{GrZu} and \cite{DiSc}.
Indeed, in those papers, the base graph $X$ is the two--way infinite
line, for which the critical percolation parameter is $\psf_c=1$, so that
the spectrum is pure point for any (finite) size of the group ``of lamps'' 
$\HH$. 
However the calculations in \cite{GrZu} and \cite{DiSc} go beyond the
present results, inasmuch the Plancherel measure is explicitly calculated.
We shall see why in general it is hard to obtain such an explicit
formula as in \cite{GrZu} and \cite{DiSc}.
Namely, the eigenfunctions of $\Rsf_{\mu}$ on $\HH \wr \GG$ arise from
eigenfunctions of each of the $P_A\,$, $A \in \Ccal$, on $\GG$.
The latter are easily computed in the case of $\GG = \Z$, since the
elements of $\Ccal$ consist of finite intervals of integers, where all
computations are explicit (compare with \cite{BaWo}). But for a general
Cayley graph $X(\GG,S)$, it appears literally impossible to compute
explicitly all eigenvalues and eigenfunctions of $P_A$ associated with
all possible finite, connected subgraphs $A$.
More explicit examples of Plancherel measures would be desirable, for example,
in order to gain more insight in the Atiyah conjecture about the possible 
values of $L^2$-Betti numbers. Besides \cite{GrZu} and \cite{DiSc}, see 
{\sc Grigorchuk, Linnell, Schick and \.Zuk}~\cite{GLSZ}, 
where this approach is explained.

Stated as above, Theorem \ref{thm:pure} can be deduced directly from Theorem
\ref{thm:spec-meas}, except for the fact that there is a complete orthonormal 
system of \emph{finitely supported} eigenfunctions. 
In physics this phenomenon is called \emph{Anderson localization}.
Below, we shall see why this is true. We shall also explain how 
the operator $\Rsf_{\wt\mu}$ can be diagonalized via the eigenfunctions of
$P_A\,$, $A \in \Ccal\,$, in the case when the percolation
clusters are a.s.\ finite.

\begin{cor}\label{cor:lattices} If $\GG = \Z^2$ and the support $S$ of the symmetric
measure $\mu$ is such that $X(\Z^2,S)$ is the square lattice or the
triangular lattice, then the associated switch-walk-switch
random walk on $\HH \wr \GG$ has pure point spectrum for any
finite group $\HH$.
\end{cor}

\begin{proof} For site percolation on the square lattice, 
$\psf_c > 1/2$. Indeed, for bond percolation, the critical probability is 
$1/2$ by a famous result of {\sc Kesten}~\cite{Kes}, and 
$\psf_c[\text{site}] > \psf_c[\text{bond}]$. 


For site percolation on the triangular lattice, $\psf_c = 1/2$.
It is also known that critical site percolation satisfies $\theta(1/2) = 0$. 

All those results can be found in \cite{Gr} and \cite{LyPe}.
\end{proof}

\begin{cor}\label{cor:tree} If $\GG$ and $S$ are such that the Cayley
graph $X(\GG,S)$ is the homogeneous tree with degree $d+1$, then
the switch-walk-switch random walk associated with any symmetric probability
measure $\mu$ supported by $S$ has pure point spectrum whenever
the group of ``lamps'' satisfies $|\HH| \ge d$.
\end{cor}

\begin{proof}
  For site percolation on the tree with degree $d+1$, 
  $\psf_c = 1/d$ and $\theta(1/d) =0$. This is an easy consequence of
  interpreting percolation in terms of a Galton-Watson process. See e.g.\ 
  \cite{LyPe}.
\end{proof}
Let us now compare the spectra of $\Rsf_{\wt\mu}$ and $\Rsf_{\mu}$ with each
other.
First we remark that since $\CC$ is amenable and $\GG= \CC \rtimes \GG / \CC$,
it follows from an old result of {\sc Kesten}~\cite[Cor. 2]{Ke} 
that the spectral radii of $\Rsf_{\wt \mu}$ and $\Rsf_{\mu}$ coincide;
see also {\sc \.Zuk}~\cite{Zu} and {\sc Woess}~\cite{Wnorms}.
More can be said if the spectrum is connected.
\begin{cor}\label{cor:spec}
  If the spectrum of the random walk on $\GG$ with law
  $\mu$ is an interval then it coincides with the spectrum of the 
  lamplighter random walk on $\HH \wr \GG$ with law $\tilde \mu$. 
  The set $\Lambda$ is dense in both spectra.
\end{cor}
\begin{proof}
The spectrum of an arbitrary selfadjoint operator $\Rsf$ on a Hilbert space 
$\mathcal H$ is contained in its \emph{numerical range}
$$
w(\Rsf) = \{ \langle \Rsf f, f \rangle : f \in \mathcal H, \norm{f}=1 \}
$$
Moreover, if the spectrum is connected, then by the min-max principle it 
coincides with the numerical range.
Recall that the spectrum of an operator does not depend on the $C^*$-algebra
which it lives in (see e.g. {\sc Takesaki}~\cite[Prop.~4.8, p.20]{Ta}). 
The base group $\GG$ is a subgroup of $\HH \wr \GG$ 
and the corresponding inclusions hold for the reduced $C^*$-algebras 
and Von Neumann algebras as well.
Thus the simple random walk transition operator $R_\mu$ can be regarded
as an element of $C_\lambda^*(\HH\wr\GG)$ and the spectrum remains the same.
Since $\Rsf_{\wt\mu}=\Rsf_\nu \Rsf_{\mu} \Rsf_\nu$ can be understood as the compression of
$\Rsf_{\mu}$ by the projection $\Rsf_\nu\,$, it is clear that 
$w(\Rsf_{\wt\mu})\subseteq w(\Rsf_{\mu})$ and therefore
$\spec(\Rsf_{\wt\mu})\subseteq \spec(\Rsf_{\mu})$.

On the other hand, $\spec \Rsf_{\mu}$ is contained in the closure of $\Lambda$.
Indeed, let $A(n) \in \Ccal$ be an increasing sequence
whose union is $\GG$. Then the $P_{A(n)}\,$, viewed as self-adjoint 
operators on $\ell^2(\GG)$, converge strongly to $\Rsf_{\mu}\,$. This implies
weak convergence of the respective resolutions of the identity, see e.g.
{\sc Dunford and Schwarz}~\cite[\S X.7]{DuSc}.
In particular, their diagonal elements at $(e,e)$, which are probability 
measures on subsets of $\Lambda$, converge weakly to the Plancherel measure of
$\Rsf_{\mu}\,$.
In conclusion, we have shown that
$\overline\Lambda \subseteq \spec \Rsf_{\wt\mu}\subseteq \spec \Rsf_{\mu}\subseteq \overline \Lambda$
and the inclusions are actually identities.
\end{proof}

A sufficient condition for the spectrum of $\Rsf_\mu$ to be connected is
the absence of nontrivial projections from the reduced $C^*$-algebra
of $\GG$. Indeed, if $I$ is a connected component of the spectrum,
the associated spectral projection can be expressed by the
analytic functional calculus, see \cite[Ch.~VII]{DuSc}
$$
\frac{1}{2\pi i} \oint_\gamma (z-\Rsf_\mu)^{-1} dz
$$
where $\gamma$ is a curve in the complement of $\spec{\Rsf_\mu}$ which
encloses $I$ and such that the
other parts of the spectrum are not inside $\gamma$.
If the spectrum is not connected, then such a spectral projection
is nontrivial.

The Kadison-Kaplansky conjecture asserts that the $C^*$-algebras
of torsion-free groups are projectionless. The Kadison-Kaplansky conjecture
in turn is implied by the surjectivity part of the Baum-Connes conjecture
\cite{Va},
which has been established for a wide range of groups,
including, e.g., amenable and Gromov hyperbolic groups.

\medskip

This paper is organized as follows: in Section~\ref{sec:proofs}, we prove the three theorems. 
We first exhibit the simple (and in principle known) argument that the
moments of the Plancherel measure of $\Rsf_{\wt \mu}$ and the integrated
density of states coincide (\S\ref{sec:proofs}.\ref{ssec:proofofthmspec-meas}). 
Then (\S\ref{sec:proofs}.\ref{ssec:projections}) we show how one can find 
projections of $\ell^2(\HH \wr \GG)$ which are sums of $\Rsf_{\wt \mu}$-invariant 
subspaces, each of which is spanned by finitely supported eigenfunctions.
This is based upon the methods of \cite{DiSc} and leads to the proofs of
theorems~\ref{thm:pointspec} 
and~\ref{thm:pure} (\S\ref{sec:proofs}.\ref{ssec:proofofthmpointspec}). 
We then proceed by explaining 
how the convolution operator can be diagonalized when the percolation clusters
are a.s.\  finite. This is quite easy when $\GG$ is torsion free (\S\ref{sec:proofs}.\ref{ssec:diagonalizationtorsionfree}),
and requires more work in the torsion case (\S\ref{sec:proofs}.\ref{ssec:diagonalizationtorsion}). 

In Section~\ref{sec:bond}, we explain how analogous results relate bond percolation
with a lamplighter random walk where the lamps are placed on the edges.

In the final Section~\ref{sec:final}, we address an open question and some possible 
extensions. In particular, we explain briefly how the (annealed, i.e., expected)
spectrum of absorbing random walk on the cluster
associated with an arbitrary percolation parameter $\psf \in (0\,,\,1)$ 
coincides with the spectrum of a deterministic convolution operator (by a signed
measure) on the wreath product $\HH \wr \GG$, where $\HH =\Z$ or another suitable
group. This is again based on the ideas of \cite{DiSc}.

In concluding the introduction, we point out that the methods and results of
this note can be nicely formulated in the more abstract language of
group $C^*$- and von Neumann algebras. Since one of our intentions is to
clarify the connection between two different circles of ideas, we have
decided to present most of the material in the more basic setup of convolution and 
random walks on groups.

We would like to thank A.~Bendikov for several stimulating discussions.
Many computations were done in the \textsc{FriCAS} computer algebra system
(a fork of the \textsc{Axiom} project, see \cite{Ax}), and we thank the members 
of the \textsc{Axiom} mailing list for their help.

\section{Proofs of the three theorems and diagonalization of the
convolution operator}\label{sec:proofs}
With $\nu$ we associate the signed measure $\nnu = \delta_{\idH} - \nu$
on $\HH$. For $g\in\GG $ we define 
$$
\nu_g = \delta_{g} * \nu * \delta_{g^{-1}} \AND 
\nnu_g = \delta_{g} * \nnu * \delta_{g^{-1}}\,.
$$
These symmetric measures live in $\HH \wr \GG\,$, and should again be 
understood in
terms of the embeddings of $\HH$ and $\GG$ into the wreath product.
It is straightforward that 
\begin{equation}\label{eq:nueqs}
\nu_g * \nu_g = \nu_g \,,\quad \nu_g * \nnu_g = 0\,,\AND
\nu_g*\nu_{g'} = \nu_{g'}*\nu_g \quad
\text{for all} \;\;g, g' \in \GG\,.
\end{equation}
For a finite set $A = \{ x_1, \dots, x_r \} \subset \GG\,$ let
\begin{equation}\label{eq:nuC}
\nu_A = \nu_{x_1} * \ldots * \nu_{x_r} \AND 
\nnu_A = \nnu_{x_1} * \ldots * \nnu_{x_r}\,.
\end{equation}
Notice that $\nnu_A\ne\delta_{(\uno,\idG)} - \nu_A$.
By the above relations, the convolution products of \eqref{eq:nuC} do not depend on the 
order (or even multiplicity) in which the elements of $A$ appear.
A short computation shows that $\nu_A$ is equidistributed on the set
$\{ (\eta,\idG) : \supp(\eta) \subset A \}$, which has $\HH^{|A|}$ elements.

\subsection{Proof of Theorem \ref{thm:spec-meas}} 
\label{ssec:proofofthmspec-meas}

Since both the Plancherel 
measure and the integrated density of states are compactly supported,
they are characterized by their moments. Thus, we only have to prove
\eqref{eq:prob-ident}. For this purpose alone, neither symmetry nor finite
support of $\mu$ are needed.

We can realize the random walk on $\GG$ with law $\mu$ on a suitable 
probability space,
and denote the corresponding probability and expectation by $\Prob$ and
$\Ex$, respectively, without adding a $\psf$ to the index.

We write $Z_n$ for the position at time $n$ of the random walk with law $\mu$ 
starting at $Z_0=\idG$. This is the product of $n$ independent $\GG$-valued random 
variables with identical distribution $\mu$. We can realize the random walk 
on a suitable probability space,
and denote the corresponding probability and expectation by $\Prob$ and
$\Ex$, respectively, without adding a $\psf$ to the index.

There is a well known and easy to prove formula for the return 
probabilities of
the switch-walk-switch random walk in terms of the random walk
on $\GG$, see e.g. {\sc Pittet and Saloff-Coste}~\cite[(2.7) and (3.1)]{PiSa}.
We recall the proof briefly.
If $g_1, \dots, g_n \in \GG$ then recursively,  
$$
\begin{gathered}
\nu * \delta_{g_1}* \nu * \delta_{g_2}* \nu * \ldots * \delta_{g_n}* \nu 
= \nu_{\idG} * \nu_{x_1} * \nu_{x_2} * \ldots * \nu_{x_n} * \delta_{x_n}
= \nu_{A} * \delta_{x_n}\,, \\\text{where}\quad 
x_k = g_1 \cdots g_k \AND A = \{ \idG, x_1, \dots, x_n \}\,.
\end{gathered}
$$
Therefore, with the same notation, and using that $\mu(g_k) = p(x_{k-1},x_k)$,
$$
\begin{aligned}
\wt\mu^{(n)} &= \nu * (\mu * \nu)^{(n)}\\
&= 
\sum_{g_1, \dots, g_n \in \GG} \Bigl(\mu(g_1)\mu(g_2)\cdots \mu(g_n)\Bigr)\,
\nu * \delta_{g_1} * \nu * \delta_{g_2} * \nu *
\ldots * \delta_{g_n} * \nu \\
&= \sum_{x_1, \dots, x_n \in \GG} 
\Bigl(p(\idG,x_1)p(x_1,x_2)\cdots p(x_{n-1},x_n)\Bigr)\,\nu_{A}*\delta_{x_n}\,.
\end{aligned}
$$
We get
$$
\begin{aligned}
\wt\mu^{(n)}(\uno,\idG) &= \sum_{x_1, \dots, x_{n-1} \in \GG}
p(\idG,x_1)p(x_1,x_2)\cdots p(x_{n-1},\idG)\, 
          \psf^{|\{ e, x_1, \dots, x_{n-1},x_n=e\}|}\\
&= \Ex\bigl(\psf^{|\{ Z_0, Z_1,\dots,Z_n \}|}\cdot \delta_{\idG}(Z_n)\bigr)\,.
\end{aligned}  
$$
(The quantity $|\{ Z_0, Z_1,\dots,Z_n \}|$, often called the \emph{range}
of the random walk,  is the number of distinct points
visited up to time $n$. In our case, $Z_n$ does not contribute to that
number, since $Z_n=Z_0$.)

Regarding the random walk on the cluster 
$C_{\omega}(\idG)$, we have that
$p_{\omega}(x_0,x_1) \cdots p_{\omega}(x_{n-1},x_n)=
p(x_0,x_1) \cdots p(x_{n-1},x_n)$ if and only if all the vertices 
$x_0, \dots, x_n$ are open, which occurs with $\Prob_{\psf}$-probability
$\psf^{|\{ x_0, x_1, \dots, x_n\}|}$. If not all of them are open then
$p_{\omega}(x_0,x_1) \cdots p_{\omega}(x_{n-1},x_n)=0\,$.
Therefore, 
$$
\Ex_{\psf}\Bigl(p_{\omega}(x_0,x_1) \cdots p_{\omega}(x_{n-1},x_n)\Bigr)
= p(x_0,x_1) \cdots p(x_{n-1},x_n) \, \psf^{|\{ x_0, x_1, \dots, x_n\}|}\,.
$$
We conclude that
$$
\Ex_{\psf}\bigl( p_{\omega}^{(n)}(\idG,\idG) \bigr)
= \Ex_{\psf} \biggl( \sum_{\,x_1, \dots, x_{n-1} \in \GG}
p_{\omega}(\idG,x_1)p_{\omega}(x_1,x_2)\cdots p_{\omega}(x_{n-1},\idG)\biggr)
$$
coincides with $\wt\mu^{(n)}(\uno,\idG)$.\qed

\smallskip

We remark that the above proof just displays the quite obvious fact that
absorbing random walk on $C_{\omega}$ is equivalent with the \emph{Rosenstock
trap model:} each point of $\GG$ is a ``trap'' with probability $(1-\psf)$,
independently of all other points. The random walk on $\GG$ survives
only before reaching a trap. In that model, one asks (e.g.) for the probability
that the random walk is back at the starting point $\idG$ in $n$ steps
without being trapped. This is just 
$$
\Ex_{\psf}\bigl( p_{\omega}^{(n)}(\idG,\idG) \bigr)
= \Ex\bigl(\psf^{|\{ Z_0, Z_1,\dots,Z_n \}|}\cdot \delta_{\idG}(Z_n)\bigr)\,.
$$
See the book by {\sc Hughes}~\cite{Hu} for many details on the
Rosenstock model, and {\sc Revelle}~\cite{Re} for its relation with
lamplighter random walks. In other words, the proof of Theorem 
\ref{thm:spec-meas} has already been implicit in the literature, but
apparently without revealing its relevance and implications within
the context addressed in the present note.
 
\subsection{Projections and invariant subspaces} 
\label{ssec:projections}
Let $A \subset \GG$ be finite. 
The (outer) \emph{vertex boundary}
$dA$ of $A$ in the Cayley graph $X(\GG,S)$ is the set of all 
$y \in \GG \setminus A$ which have a neighbour in $A$. When $A = \emptyset$,
we define $dA = \{ \idG \}$. We introduce the 
signed measures
\begin{equation}\label{eq:nuCdA}
\nu_{A,dA} = \nu_A * \nnu_{dA}\,.
\end{equation}

We write $Q_A$ for the right convolution operator by $\nu_{A,dA}$.
It is a projection of $\ell^2(\HH \wr \GG)\,$ in that 
$\nu_{A,dA} * \nu_{A,dA}=\nu_{A,dA}$ by  \eqref{eq:nueqs}. We write 
$$
Q_A\ell^2 = \{ F * \nu_{A,dA} : F \in \ell^2(\HH \wr \GG) \}\,.
$$
for its image. 

\begin{pro}\label{pro:decompose} \emph{(a)} The projections $Q_A\,$, 
$A \in \Ccal$,
are mutually orthogonal.\\[3pt]
\emph{(b)} We have 
$$
\ell^2(\HH \wr \GG) = \bigoplus_{A \in \Ccal} Q_A\ell^2
\quad\; \text{(closed direct sum), or equivalently,} \;\quad
\delta_{(\uno, \idG)} = \sum_{A \in \Ccal} \nu_{A,dA}
$$
if and only if $\theta(\psf)=0$ for site percolation on $X(\GG,S)$ with 
parameter $\psf=1/|\HH|\,$.
\end{pro}

\begin{proof} (a) If $A, B \in \Ccal$ and $A \ne B$ then
both $A$ and $B$ are connected and have $e$ as common vertex.
Any path from $e$ to a vertex in say $A\setminus B$ must cross $dB$.
It follows that at least one of $A \cap dB$ and $dA \cap B$ is non-empty,
Let $x$ be an element of one of those sets. Then, using \eqref{eq:nueqs}, we see that
$\nu_{A,dA} * \nu_{B,dB}$ contains the convolution factor $\nu_x * \nnu_x = 0$.
Therefore $Q_BQ_A = 0$.   

(b) It follows from the general theory 
of von Neumann traces (see e.g. {\sc Takesaki}~\cite{Ta}), 
here specifically in
the context of group algebras, that the statement
is equivalent with
\begin{equation}\label{eq:trace-sum}
\sum_{A \in \Ccal} \nu_{A,dA}(\uno,\idG) = 1\,.
\end{equation}
As a matter of fact, it does not require full understanding of
von Neumann traces to see why \eqref{eq:trace-sum} suffices for the
proof. Let $\phi = \sum_{A \in \Ccal} \nu_{A,dA}$
and $\psi = \delta_{(\uno, \idG)} - \phi$. Then $\phi$ and $\psi$ are
orthogonal projections, and in particular, they are both positive definite.
This implies in particular that the real, symmetric matrix 
$
\left(\begin{smallmatrix} \psi(\uno, \idG) & \psi(\eta, g) \\
                     \psi(\eta, g) & \psi(\uno, \idG)
\end{smallmatrix}\right)
$ 
a compression of $\phi$ and thus positive definite for each $(\eta, g) \in \HH \wr \GG\,$
different from $(\eta,g)$,
that is
$\psi(\uno, \idG)^2 \ge \psi(\eta, g)^2\,.$ 
Thus $\psi \equiv 0$ if and only
if $\phi(\uno, \idG) = 1\,.$ 		   

We come back to the verification of \eqref{eq:trace-sum}: since 
$$
\nu_{A,dA}(\uno,\idG) = \psf^{|A|} (1-\psf)^{|dA|} 
= \Prob_{\psf}[C(\idG) = A]\,,
$$
the proposition follows.
\end{proof}

\begin{lem}\label{lem:invariant} Let $A \in \Ccal$ be non-empty.  
Then for $f: \GG \to \R$ with  $\supp(f) \subset A$, 
$$
\nu_{A,dA} * f * \wt \mu = \nu_{A,dA} * P_Af\,.
$$
Thus, the linear span 
$$
\Bigl\{ \nu_{A,dA} * f \;:\; f \in \R^{\GG}\,,\; \supp(f) \subset A \Bigr\}
$$
is mapped into itself under right convolution with $\wt \mu$. 
\end{lem}
 
\begin{proof}
For $x \in A$ and $s \in S$
$$
\nu_{A,dA} * \delta_x * \nu * \delta_s * \nu 
= \nu_A * \nnu_{dA} * \nu_x * \nu_{xs} * \delta_{xs} =
\begin{cases} \nu_A * \nnu_{dA}* \delta_{xs}\,,&\text{if}\; xs \in A\,,\\
              0\,,&\text{if}\; xs \in dA\,.
\end{cases}
$$	      
Therefore, for $x \in A$,
$$
\begin{aligned}
\nu_{A,dA} * \delta_x * \wt\mu 
&= \sum_{s \in S} \mu(s) \cdot \nu_{A,dA} * \delta_x * \nu * \delta_s * \nu \\
&= \nu_{A,dA} *  
\biggl(\sum_{\,s \in S: xs \in A} \mu(s)\cdot \delta_{xs}\biggr)
= \nu_{A,dA} * (P_A \delta_x)\,.
\end{aligned}
$$
(We have used symmetry of $\mu$ in the last step.)
The result now follows simply by writing $f=\sum_x f(x)\cdot\delta_x$.
\end{proof}

We see that for non-empty $A \in \Ccal$, the space
$$
\wt Q_A \ell_2 = 
\bigl\{ Q_A(F) * f : F \in \ell^2(\HH \wr \GG)\,,\;
f \in \R^{\GG}\,,\; \supp(f) \subset A \bigr\}
$$
is mapped into itself by the right convolution operator $\Rsf_{\wt \mu}$.

In the specific case when $A = \emptyset$ then 
$dA = \{ \idG \}$ by definition, and $\nu_{A,dA} = \nnu_{\idG}$. 
We have $\nnu_{\idG} * \wt\mu = 0\,$, so that $Q_{\emptyset}\ell^2$
is mapped to $\{ 0\}$ by the right convolution operator $\Rsf_{\wt \mu}$.

Now, again for non-empty $A$, let $\{ f_{A,x} : x \in A\}$ be an orthonormal
system of right eigenfunctions of the symmetric matrix $P_A$ with 
associated  eigenvalues $\la_{A,x} \in \R$. We write 
\begin{equation}\label{eq:SAx}
\sigma_{A,x} = \nu_{A,dA} * f_{A,x}\,, \AND S_{A,x}(F) = F * \sigma_{A,x}\,, 
\end{equation}
$F \in \ell^2(\HH \wr \GG)$, for the associated right convolution operator.
Then by Lemma \ref{lem:invariant}
\begin{equation}\label{eq:sigma}
\sigma_{A,x} * \wt \mu 
= \la_{A,x} \cdot \sigma_{A,x}\,.
\end{equation}

\subsection{Proof of theorems \ref{thm:pointspec} and \ref{thm:pure}}
\label{ssec:proofofthmpointspec}
Let $A \in \CC$ be non-empty. Then $\wt Q_A \ell^2$ is the sum of its finitely
many subspaces
$$
S_{A,x}\ell^2 =  
\bigl\{ \underbrace{Q_A(F)*f_{A,x}}_{\displaystyle S_{A,x}(F)} : 
   F \in \ell^2(\HH \wr \GG)\bigr\}\,,\; x \in A\,.
$$
By \eqref{eq:sigma}, $S_{A,x}\ell^2$ is an eigenspace of $\Rsf_{\wt \mu}$ 
with eigenvalue $\la_{A,x}$. It is generated by all functions
$S_{A,x}(F)$, where $\supp(F)$ is finite. Since $S_{A,x}$ is a convolution
operator by a finitely supported signed measure, all those functions are 
finitely supported, and $S_{A,x}\ell^2$ is infinite-dimensional.

Furthermore, we have the eigenspace $Q_{\emptyset}\ell^2$ with eigenvalue
$0$. It is infinite-dimensional and generated by finitely supported
functions by the same reason as above. For convenience, we write
$\wt Q_{\emptyset}\ell^2 = Q_{\emptyset}\ell^2\,$.

This proves Theorem \ref{thm:pointspec}. 

Next, suppose that $\theta(\psf)=0$ for site percolation with parameter 
$\psf = 1/|\HH|\,$.
Note that for arbitrary $y \in A$ we have
$\nu_{A,dA} * \delta_y = \delta_y * \nu_{B,dB}\,$, where 
$B = y^{-1}A$ is again in $\Ccal\,$. This identity implies
that 
$$
Q_A \ell^2 * \delta_y = \{ F * \delta_y : F \in Q_A \ell^2 \}
= Q_{y^{-1}A} \ell_2\,.
$$
Since every $\delta_y$, $y \in A$, can be written as a linear combination
of the functions $f_{A,x}\,$, $x \in A$, we see that 
$\wt Q_A$ is also the sum of its subspaces $Q_{y^{-1}A} \ell_2\,$, 
$y \in A$. As $A \in \Ccal$ and (if $A \ne \emptyset$) the $y \in A$ vary, 
it follows from Proposition~\ref{pro:decompose} that
the sum of the latter spaces is dense in $\ell^2(\HH \wr \GG)$. 
We conclude that the same is true for the union of the spaces $\wt Q_A\,.$
By the above, the latter are generated by finitely supported
eigenfunctions of the convolution operator $\Rsf_{\wt \mu}\,$. This
concludes the proof of Theorem \ref{thm:pure}. \qed

\subsection{Diagonalization of the convolution operator in the torsion-free
case}
\label{ssec:diagonalizationtorsionfree}

We next want to describe a diagonalization of the operator $\Rsf_{\wt \mu}$
associated with the switch-walk-switch random walk.

Recall the operator $S_{A,x}$ of \eqref{eq:SAx}, where $A \in \Ccal$ is
non-empty and $x \in A$. Since the measure $\nu_{A,dA}$ is symmetric,
the adjoint $S^*_{A,x}$ is the right convolution operator by 
$\check f_{A,x} * \nu_{A,dA}\,$, where $\check f(g) = f(g^{-1})$ for
$f \in \R^{\GG}\,$.

\begin{pro}\label{pro:ort} Let $A \in \Ccal$ be non-empty and such that
$g_1^{-1}A \ne g_2^{-1}A$ for all distinct $g_1, g_2 \in A$.
Then for $x, y \in A$
$$
S^*_{A,x}\, S_{A,y} = \begin{cases} Q_A\,,&\text{if}\; x=y\\
                                  0\,,&\text{otherwise.}
                    \end{cases}
$$
\end{pro}

\begin{proof} 
Let  $g_1, g_2 \in A$. Set $B = g_1^{-1}A$ and
$C=g_2^{-1}A$. We use \eqref{eq:nueqs} and Proposition \ref{pro:decompose}(a)
to compute
$$
\nu_{A,dA} * \delta_{g_1} * \delta_{g_2^{-1}} * \nu_{A,dA}
= \delta_{g_1} * \nu_{B,dB} * \nu_{C,dC} * \delta_{g_2^{-1}} =
 \begin{cases} \delta_{g_1} * \nu_{B,dB} *  \delta_{g_2^{-1}} 
                        \,,&\text{if}\; B=C     \\
                                  0\,,&\text{otherwise.}
                    \end{cases}
$$
Now $B=C$ means that $g_1^{-1}A = g_2^{-1}A$, which by assumption implies
$g_1=g_2$. But then 
$\delta_{g_1} * \nu_{B,dB} *  \delta_{g_1^{-1}} = \nu_{A,dA}\,$. Therefore
the measure that induces the right convolution operator $S^*_{A,x}\, S_{A,y}$ is
$$
\begin{aligned}
\nu_{A,dA} * f_{A,y} * \check f_{A,x} * \nu_{A,dA} 
&= \sum_{g_1,g_2 \in A} f_{A,y}(g_1)f_{A,x}(g_2) \,\cdot\, 
         \nu_{A,dA} * \delta_{g_1} * \delta_{g_2^{-1}} * \nu_{A,dA}\\
&= \sum_{g \in A} f_{A,x}(g)f_{A,y}(g) \,\cdot \,\nu_{A,dA} 
= \begin{cases} \nu_{A,dA}\,,&\text{if}\; x=y\\
                                  0\,,&\text{otherwise}
                    \end{cases}
\end{aligned}
$$
by orthonormality.
\end{proof}

Now note that $g_1^{-1}A = g_2^{-1}A$ for distinct $g_1, g_2 \in A$
implies that $g_1g_2^{-1}$ stabilizes the finite set $A$ and must be a torsion
element of $\GG$. 

We now suppose for the rest of this sub-section that $\GG$ is
torsion-free.
Then the conclusion of Proposition \ref{pro:ort} is always valid.
In particular, 
$$
S_{A,x}\, S^*_{A,x}\, S_{A,x} = S_{A,x}\, Q_A = S_{A,x}\,.
$$
That is, $S_{A_x}$ is a \emph{partial isometry,} and
$$
T_{A,x} = S_{A,x}\, S^*_{A,x}\,,\quad 
T_{A,x}(F) = F * \check f_{A,x} * 
\underbrace{\nu_{A,dA} * f_{A,x}}_{\displaystyle \sigma_{A,x}}
\;\; \text{for}\;\; F \in \ell^2(\HH \wr \GG)\,,
$$is the orthogonal projection of $\ell^2(\HH \wr \GG)$ onto the (closed) 
subspace $S_{A,x}\ell^2\,$. We subsume.

\begin{thm}\label{thm:diag} Let $\GG$ be torsion-free.
Then all the operators $T_{A,x}\,$, where 
$A$ varies in $\Ccal \setminus \{ \emptyset \}$ and $x$ varies in $A$,
together with $Q_{\emptyset}$ are mutually orthogonal projections
of $\ell^2(\HH \wr \GG)$ onto eigenspaces of the convolution operator
$\Rsf_{\wt \mu}\,$. 

All those eigenspaces are infinite-dimensional and
spanned by finitely supported functions. 
Regarding the associated eigenvalues, we have 
$$
\Rsf_{\wt \mu}\,T_{A,x} = \la_{A,x}\cdot T_{A,x} \AND \Rsf_{\wt \mu}\,Q_{\emptyset} =
0\,.
$$

If $\theta(\psf)=0$ for bond percolation on $X(\GG,S)$ with parameter 
$\psf = 1/|\HH|$ then the (closed) direct sum of those eigenspaces is the
whole of $\ell^2(\HH \wr \GG)$, providing a complete diagonalization of
the operator $\Rsf_{\wt \mu}\,$.
\end{thm}

\subsection{Diagonalization in the torsion case}

\label{ssec:diagonalizationtorsion}
If $\GG$ has torsion elements then the proof of Proposition \ref{pro:ort}
does not work anymore, and we need a more refined method for diagonalizing
$\Rsf_{\wt \mu}\,$. A somewhat similar construction for Gelfand pairs is
done in {\sc Scarabotti and Tolli}~\cite{ScTo2}.
We reformulate the diagonalization procedure for the torsion-free case
in terms of partial isometries.
Recall that we consider right convolution operators
$$
\Rsf_f \, g(x) = g *\refl f(x) = \sum_y g(xy)f(y)
$$
where $f:G\to \IC$ and $\refl{f}$ is the reflection $\refl f(x) = f(x^{-1})$.
With this convention we have
$
\Rsf_f \Rsf_g = \Rsf_{f*g}\,;
$
in other words, $R$ is the right regular representation.
For group elements $x$ we will abbreviate 
$\Rsf_x=\Rsf_{\delta_x}\,$, which is the right convolution with $\delta_{x^{-1}}$.
Right convolution with the idempotent measure $\nu$ leads to a projection 
$E=\Rsf_\nu\,$, and the translated projections $E_g = \Rsf_{\nu_g}$ satisfy 
the fundamental commutation relation $E_g = \Rsf_g\, E\, \Rsf_{g^{-1}}\,$.
Thus the projections $Q_A = \Rsf_{\nu_{A,dA}}$ can be written as
$$
Q_A = \prod_{x\in A} E_x \prod_{y\in dA} (I-E_y)
$$
and the relation 
$
\Rsf_x \,Q_A\, \Rsf_{x^{-1}} = Q_{xA}
$
holds.
Then the essence of Lemma~\ref{lem:invariant} can be recaptured algebraically as follows:
For  $A\in \Ccal$ the space spanned by the partial isometries 
$\{\Rsf_{x^{-1}}Q_A : x \in A\}$ is invariant under $\Rsf_{\wt\mu}$,
namely for $x\in A$
$$
\Rsf_{\wt \mu} \Rsf_{x^{-1}} Q_A = \sum_{y\in A} p(x,y) \Rsf_{y^{-1}}Q_A\,.
$$
In other words, if $v=(v_a)_{a\in A}$ is a vector and $W_v=\sum_{a\in A} v_a\Rsf_{a^{-1}}Q_A$,
then
\begin{equation}\label{eq:RmuWv}
  \Rsf_{\wt\mu} W_v = W_{P_Av}\,.
\end{equation}
Diagonalizing $P_A\,$, which is symmetric, yields a unitary matrix
$U = ( u_{x,y})_{x,y\in A}$ such that
$U^{-1}P_AU=\diag(\lambda_x)_{x\in A}$, i.e., the columns of $U$
are eigenfunctions of $P_A$ and the ranges of the operators
$$
S_{A,x} = \sum_y u_{y,x}\,\Rsf_{y^{-1}}\,Q_A
$$
are eigenspaces of $\Rsf_{\wt\mu}$.
In the torsion-free case the stabilizer $\GA=\{x\in \GG: xA=A\}$ is trivial
and each $S_{A,x}$ is a partial isometry whose range
projection $S_{A,x}\,S_{A,x}^*$ is a spectral projection of $\Rsf_{\wt\mu}$.
If $\GG$ is not torsion free then $\GA$ can be nontrivial and and $S_{A,x}$
need not be a partial isometry anymore. In this case finite Fourier
analysis is needed.

We recall a few facts from finite noncommutative harmonic analysis;
for details see e.g. {\sc Serre}~\cite{Se}.
Let $\Gamma$ be a finite group.
Every finite dimensional unitary representation $\rho:\Gamma\to \GL(V_\rho)$ can be decomposed into a direct sum
$$
V_\rho = \bigoplus_{\pi\in \widehat\Gamma}  m_\pi V_\pi
$$
where $\widehat\Gamma$ denotes the set of irreducible unitary representations of $\Gamma$ and
$m_\pi$ denotes the multiplicity of $\pi:\Gamma\to \GL(V_\pi)$.
In particular, the left regular representation decomposes as
$$
\lambda = \bigoplus_{\pi\in \widehat\Gamma}  d_\pi \pi
$$
where $d_\pi=\dim V_\pi$ denotes the dimension of the irreducible representation $\pi$.
For the group ring $\IC \Gamma = \linspan \lambda(\Gamma)$,
this means that there is a decomposition into a direct sum of matrix algebras
$$
\IC\Gamma \simeq \bigoplus_{\pi\in \widehat\Gamma} \linspan \pi(\Gamma)
$$
where by irreducibility each $\linspan \pi(\Gamma)$ is isomorphic to the full matrix algebra $M_{d_\pi}(\IC)$.
This isomorphism is implemented by the \emph{finite Fourier transform},
namely for a function $f:\Gamma\to\IC$ we denote for $\pi\in \widehat\Gamma$ the 
Fourier transform
$$
\hat f(\pi) = \sum_{x\in\Gamma} f(x) \, \pi(x) \in M_{d_\pi}(\IC)
$$
and the inverse Fourier transform is given by the identity
$$
f(x) = \frac{1}{\abs{\Gamma}} \sum_{\pi\in\widehat\Gamma} d_\pi 
     \Tr\bigl(\pi(x)^*\hat{f}(\pi)\bigr)\,.
$$
Writing out the latter formula, if we denote by $\pi_{ij}(x)$ the matrix 
entries of $\pi(x)$, we have
$$
f(x) = \frac{1}{\abs{\Gamma}} \sum_{\pi\in\widehat\Gamma} 
\sum_{i,j =1}^{d_\pi} d_\pi \,\overline{\pi_{ij}(x)} \, \hat{f}(\pi)_{ij}\,.
$$
In particular, if we equip $\IC G$ with the standard scalar product
$
\langle f,g \rangle = \sum_{x\in\Gamma} f(x)\overline{g(x)}\,
$,
then the matrix coefficients
\begin{equation}
  \label{eq:Fourierbasis}
  \{
  e_{\pi,i,j} = \sqrt{\tfrac{d_\pi}{\abs{\Gamma}}} \,\pi_{ij} 
  :
  \pi\in\widehat\Gamma, 1\leq i,j\leq d_\pi
  \}
\end{equation}
form an orthonormal basis and the Fourier transform implements the unitary change
from the canonical basis $\{\delta_x : x\in \IC \Gamma\}$ to the Fourier basis.
The following proposition summarizes the basic facts about finite Fourier
analysis which are needed in the sequel.
\begin{pro}
  \label{prop:schurlemma}
  Let $\pi,\rho\in\widehat\Gamma$ and
  $1\leq i,j\leq d_\pi$,
  $1\leq s,t\leq d_\rho$,
  $x,y\in\Gamma$.
  Then
  \begin{subequations}
    \label{eq:representation}
    \begin{align}
      \pi_{ij}(xy)     &= \sum_{k=1}^{d_\pi} \pi_{ik}(x)\,\pi_{kj}(y)     \\
      \pi_{ij}(x^{-1}) &= \overline{\pi_{ji}(x)}
    \end{align}
  \end{subequations}
  and by \emph{Schur's Lemma}
  \begin{subequations}
    \label{eq:schurlemma}
    \begin{align}
      \pi_{ij}*\rho_{st} &= 0      & \text{if $\pi\ne\rho$}\\
      \pi_{ij}*\pi_{st}  &= \delta_{js}\, \tfrac{\abs{\Gamma}}{d_\pi}\, \pi_{it}
    \end{align}
  \end{subequations}
  Moreover,
  \begin{equation}
    \label{eq:allrepresentations}
    \sum_{\pi\in\widehat\Gamma} \sum_{i,j=1}^{d_\pi} \pi_{ij} = \delta_e\,.
  \end{equation}
\end{pro}
In other words, the functions $e_{\pi,i,j}$ form an orthogonal family of
matrix units.

Assume now that the stabilizer $\GA=\{x\in G: xA=A\}$ is nontrivial.
$\GA$ is finite and its left action commutes with the action of 
$\Rsf_{\wt\mu}\,$, since $P_A$ is $\GA$-invariant: if $gA=A$, then 
$p_A(gx,gy)=p(gx,gy)=p(x,y)=p_A(x,y)$, that is,
$
\langle P_A \delta_{gy},\delta_{gx} \rangle 
= \langle P_A \delta_y,\delta_x\rangle.
$
We can decompose $A$ into a  finite set of disjoint orbits
$$
A = \bigcup_{k=1}^m \GA\, a_k
$$
where $\{a_k:1\leq k\leq m\}$ is a set of representatives.
Let us compute the entries of $P_A$ in the orbit-wise Fourier basis  
$$
e_{k,\pi,i,s} = \sqrt{\tfrac{d_\pi}{\abs{\GA}}} \,\pi_{is}*\delta_{a_k}
$$
where $1\leq k\leq m$, $\pi\in\hatGA$ and $1\leq i,s\leq d_\pi$.
Since the support of $e_{k,\pi,i,s}$ is the orbit $\GA a_k$, we obtain
\begin{align*}
  \langle P_A e_{k,\pi,i,s} \,,\,e_{l,\rho,j,t}\rangle
  &= \frac{\sqrt{d_\pi d_\rho}}{\abs{\GA}}
     \sum_{x,y} \pi_{is}(x)\, \overline{\rho_{jt}(y)} \,
       \langle P_A \delta_{xa_k},\delta_{ya_l} \rangle \\
  &= \frac{\sqrt{d_\pi d_\rho}}{\abs{\GA}}
     \sum_{z} \rho_{tj}* \pi_{is}(z) \,\langle P_A \delta_{za_k},\delta_{a_l} \rangle \\
  &= \delta_{\pi\rho}\, \delta_{ij} \sum_{z\in \GA} \pi_{ts}(z)\, p(a_l,za_k)
\end{align*}
and the matrix coefficients do not depend on $i$.
Thus for fixed $\pi,i$ the linear subspace 
$\linspan\{e_{k,\pi,i,s} : k=1,\dots, m\,;\,\, s=1,\dots, d_\pi\}$
is invariant and 
$$
P_A e_{k,\pi,i,s} 
= \sum_{l=1}^m \sum_{t=1}^{d_\pi} m^{(\pi)}_{lt,ks}\,
e_{l,\pi,i,t}\,,\quad\text{where}\quad
m^{(\pi)}_{ks,lt} = \sum_{z\in \GA} \pi_{ts}(z)\, p(a_l,za_k)\,,
$$
form a hermitian matrix $M^{(\pi)}$.
The basis change is given by the Fourier matrix
$$
F = \Bigl( e_{k,\pi,i,s}(xa_l) \Bigr)_{\substack{1\leq l\leq m,\, x\in \GA \\
1\leq k\leq m,\, \pi\in\hatGA,\, 1\leq i,s\leq d_\pi}}
$$
and $P_A=FMF^*$.
Since the matrix $M^{(\pi)}$ does not depend on the index $i$,
we can diagonalize $M=\bigoplus M^{(\pi)}$ by a block-diagonal unitary matrix
$U=\bigoplus U^{(\pi)}$ whose entries
$$
\bigl\{u_{k,\pi,i,s ; \, l,\rho,j,t}=\delta_{\pi\rho}\, 
\delta_{ij} \, u^{(\pi)}_{ks,lt} :
 k,l =1,\dots, m,\,\, \pi,\rho \in \hatGA\,,\,\, i,s = 1,\dots, d_\pi\,,\,\, 
 j, t = 1,\dots, d_\rho\bigr\}
$$
do not depend on $i$ either.
Then we have $P_A = FU\Lambda U^*F^*\,$, and the eigenvectors of $P_A$ are
$v_{k,\pi,i,s}=FU\delta_{k,\pi,i,s}$. By \eqref{eq:RmuWv}, the range of
$W_{v_{k,\pi,i,s}}$ is an eigenspace of $\Rsf_{\wt\mu}\,$,
and since the eigenvalues do not depend on the index $i$, the range of
$\sum_{i=1}^{d_\pi} W_{v_{k,\pi,i,s}}$ is invariant as well.
\begin{pro}\label{pro:dec}
The operators 
$$
S_{k,\pi,s} = \frac{\sqrt{d_\pi}}{\abs{\GA}} \,
\sum_{l=1}^m \sum_{i,t=1}^{d_\pi} u_{lt,ks}^{(\pi)} \, \Rsf_{a_l^{-1}} \,
\Rsf_{\overline{\pi_{ti}}}\, Q_A
$$
with $k \in \{1, \dots, m\}$, $\pi\in \hatGA$ and $s \in \{1, \dots, d_\pi\}$
form a family of partial isometries whose range projections 
$$
T_{k,\pi,s}=S_{k,\pi,s} \, S_{k,\pi,s}^*
$$
are mutually orthogonal eigenprojections of $\Rsf_{\wt\mu}$ and 
$$
\sum_{k,\pi,s} T_{k,\pi,s} = \wt{Q}_A\,.
$$
\end{pro}

\begin{proof}
Indeed $S_{k,\pi,s}=\frac{1}{\sqrt{d_\pi}} \sum_{i=1}^{d_\pi} W_{v_{k,\pi,i,s}}$
with 
\begin{align*}
    v_{k,\pi,i,s}(xa_p) &= FU\delta_{k,\pi,i,s} \\
                 &= \sum_{l=1}^m \sum_{\rho\in\hatGA} \sum_{j,t=1}^{d_\rho} 
                     F_{p,x;\,l,\rho,j,t} \,U_{l,\rho,j,t;\,k,\pi,i,s} \\
                 &= \sum_{l=1}^m \sum_{\rho\in\hatGA} \sum_{j,t=1}^{d_\rho} 
                     e_{l,\rho,j,t}(xa_p) \, \delta_{\pi\rho}\, 
		             \delta_{ij} \, u^{(\pi)}_{lt,ks} \\
                 &= \sum_{l=1}^m  \sum_{t=1}^{d_\pi} 
                     \delta_{pl} \sqrt{\tfrac{d_\pi}{\abs{\GA}}} 
		           \,\pi_{it}(x) \, u^{(\pi)}_{lt,ks} \\
                 &= \sum_{t=1}^{d_\pi} 
                     \sqrt{\tfrac{d_\pi}{\abs{\GA}}}\, \pi_{it}(x) 
		            \, u^{(\pi)}_{lt,ks} 
\end{align*}
By the previous calculations it follows that 
$\Rsf_{\wt\mu} \, S_{k,\pi,s}= \lambda_{ks}^{(\pi)}\cdot S_{k,\pi,s}$ 
and it remains to check  orthogonality and the partial isometry condition:
  \begin{align*}
    S_{k,\pi,s}^*\, S_{k',\pi',s'} 
    &= \frac{\sqrt{d_\pi d_{\pi'}}}{\abs{\GA}^2}
       \sum_{l=1}^m \sum_{i,t=1}^{d_\pi} 
       \sum_{l'=1}^m \sum_{i',t'=1}^{d_{\pi'}}
       \overline{u_{lt,ks}^{(\pi)}}\,
       u_{l't',k's'}^{(\pi')}\,
       Q_A \,\Rsf_{\overline {\pi_{ti}}}^*\, \Rsf_{a_l}\,
       \Rsf_{a_{l'}^{-1}}\, \Rsf_{\overline {\pi'_{t'i'}}}\,Q_A\\
\intertext{noting that $\Rsf_{\overline{\pi_{ti}}}$ commutes with $Q_A$ and 
$Q_A\Rsf_{a_l} \,\Rsf_{a_{l'}^{-1}}\,Q_A=\delta_{ll'}\,Q_A$ we get}
    &= \frac{\sqrt{d_\pi d_{\pi'}}}{\abs{\GA}^2}
       \sum_{l=1}^m \sum_{i,t=1}^{d_\pi} 
       \sum_{i',t'=1}^{d_{\pi'}}
       u_{ks,lt}^{(\pi)}\,
       u_{l't',k's'}^{(\pi')}\,
       \Rsf_{\overline {\pi_{it}*\pi'_{t'i'}}}\,Q_A\\
\intertext{using orthogonality of the Fourier basis
\eqref{eq:schurlemma}
this simplifies to
}
    &= \frac{1}{\abs{\GA}}
       \sum_{l=1}^m \sum_{t,i,i'=1}^{d_\pi} 
       u_{ks,lt}^{(\pi)}\,
       u_{lt,k's'}^{(\pi)}\,
       \Rsf_{\overline {\pi_{ii'}}}\, Q_A
\intertext{and by the unitary property 
$\sum_{l=1}^m \sum_{t=1}^{d_\pi} u_{ks,lt}^{(\pi)}\, u_{lt,k's'}^{(\pi)} 
= \delta_{kk'}\,\delta_{ss'}$ we have}
    &= \delta_{\pi\pi'}\, \delta_{kk'} \,\delta_{ss'}\, \frac{1}{\abs{\GA}}
       \sum_{i,i'=1}^{d_\pi} 
       \Rsf_{\overline {\pi_{ii'}}}\, Q_A
  \end{align*}
  and the result is indeed a projection.
  Next we check the partial isometry condition:
  \begin{align*}
    S_{k,\pi,s} \, S_{k,\pi,s}^* \, S_{k,\pi,s}
    &= \frac{\sqrt{d_\pi}}{\abs{\GA}^2}\,
       \sum_{l=1}^m \sum_{i,t,j,j'=1}^{d_\pi} 
       u_{lt,ks}^{(\pi)}\,
       \Rsf_{a_{l}^{-1}} \,\Rsf_{\overline {\pi_{ti}}}\,
       \Rsf_{\overline {\pi_{jj'}}}\,Q_A \\
    &= \frac{\sqrt{d_\pi}}{\abs{\GA}^2}\,
       \sum_{l=1}^m \sum_{i,t,j,j'=1}^{d_\pi} 
       u_{lt,ks}^{(\pi)}\,
       \Rsf_{a_{l}^{-1}} \,\delta_{ij}\, \frac{\abs{\GA}}{d_\pi} \,
          \Rsf_{\overline {\pi_{tj'}}}\, Q_A\\
    &= \frac{\sqrt{d_\pi}}{\abs{\GA}}\,
       \sum_{l=1}^m \sum_{t,j'=1}^{d_\pi} 
       u_{lt,ks}^{(\pi)}\,
       \Rsf_{a_{l}^{-1}}\, \Rsf_{\overline {\pi_{tj'}}}\, Q_A\\
    &= S_{k,\pi,s}\,.
  \end{align*}
  Finally, the sum of the range projections is
  \begin{align*}
  \sum_{k=1}^m &\sum_{\pi\in\hatGA} \sum_{s=1}^{d_\pi} 
  S_{k,\pi,s} \,S_{k,\pi,s}^*\\
  &= \sum_{\pi\in\hatGA} \sum_{k=1}^m \sum_{s=1}^{d_\pi}
     \frac{d_\pi}{\abs{\GA}^2}
       \sum_{l=1}^m \sum_{i,t=1}^{d_\pi} 
       \sum_{l'=1}^m \sum_{i',t'=1}^{d_\pi} 
       u_{lt,ks}^{(\pi)}\,
       \overline{u_{l't',ks}^{(\pi)}}\,
       \Rsf_{a_{l}^{-1}} \,\Rsf_{\overline {\pi_{ti}}}\, Q_A\,
       \Rsf_{\overline {\pi_{t'i'}}}^*\, \Rsf_{a_{l'}}\\
  &= \sum_{\pi\in\hatGA} 
     \frac{d_\pi}{\abs{\GA}^2}
       \sum_{l=1}^m \sum_{i,t=1}^{d_\pi} 
       \sum_{l'=1}^m \sum_{i',t'=1}^{d_\pi} 
       \underbrace{
         \sum_{k=1}^m \sum_{s=1}^{d_\pi}
         u_{lt,ks}^{(\pi)}\,
         u_{ks,l't'}^{(\pi)}
       }_{\delta_{ll'}\delta_{tt'}}\,
       \Rsf_{a_{l}^{-1}} \,
       \delta_{ii'}\, \frac{\abs{\GA}}{d_\pi}\, \Rsf_{\overline {\pi_{tt'}}}\, 
       Q_A\,
       \Rsf_{a_{l'}}\\
  &= \frac{d_\pi}{\abs{\GA}}
       \sum_{l=1}^m 
       \Rsf_{a_{l}^{-1}} 
       \sum_{\pi\in\hatGA} 
       \sum_{t=1}^{d_\pi} 
       \Rsf_{\overline {\pi_{tt'}}} \,
       Q_A\,
       \Rsf_{a_{l'}}
  = \sum_{l=1}^m 
       Q_{a_l^{-1}A} 
  = \wt Q_A\,,
  \end{align*}
  because of 
  \eqref{eq:allrepresentations}.
\end{proof}

\section{Bond percolation and enlightened edges}\label{sec:bond}

At this point, it appears natural to ask whether there is an analogous
relation between \emph{bond} percolation and a suitable
lamplighter random walk. 

\subsection{Bond percolation}
We need the (unoriented) edge set 
$\Ef = \{ [x,xs] = [xs,x] : x \in \GG\,,\; s \in S \}$ of $X(\GG,S)$.
We do not repeat all basic details: bond percolation is analogous to
what is explained in \S\ref{sec:intro}.\ref{ssec:percolationclusters}. 
The only difference is that we now have i.i.d.
Bernoulli random variables $Y_{\ef}$, where $\ef \in \Ef\,$. An edge $\ef$ 
is open if   
$Y_{\ef} = 1$, which occurs with probability $\psf$, and closed, otherwise. 
The closed edges are
removed, and the random clusters are the connected components of the
graph that is left over. 

One slight difference arises in the terminology of connected \emph{subgraphs}.
In the preceding sections, they were \emph{induced subgraphs,} where 
we have $A \subset \GG$ and turn it into a subgraph of
$X(\GG,S)$ by keeping all edges of the latter which have both
endpoints in $A$. In the present context, connected subgraphs are not 
necessarily induced; they are just such that their vertex and edge sets
are subsets of $\GG$ and $\Ef$. For such a subgraph $A$, we write $\Ef(A)$ for its
set of edges, and $\partial A$ for the set of edges in $\Ef$ that do not belong 
to $\Ef(A)$ but have almost surely an endpoint in $A$. The restricted 
transition probabilities of the random walk on $\GG$ with law $\mu$ now become
\begin{equation}\label{eq:restrRW}
p_A(x,y) = \begin{cases} p(x,y)=\mu(x^{-1}y)\,,&\text{if}\; 
                            [x, y] \in \Ef(A)\,,\\
                       0\,,&\text{otherwise.}
         \end{cases}
\end{equation}
We have to replace $\Ccal$
with the family $\Ccal^{\ed}$ of all finite connected subgraphs (in this wider 
sense) of $X(\GG,S)$ containing $\idG$. If  
$A= C_{\omega}^{\ed}(\idG)$ is the cluster of $\idG$ then we write
$p_{\omega}(x,y) = p_{A}(x,y)$. If $A$ is finite then
it belongs to $\Ccal^{\ed}\,$. 
Here, we do not need to add the empty set to $\Ccal^{\ed}\,$;
if percolation is such that all edges incident with $\idG$ are closed, then 
the resulting component is $\{ \idG \}$. Elements of $\Ccal^{\ed}$ are called
\emph{animals} in \cite{Gr}.

Taking into account those modifications,
the definitions of the spectral measures $\nsf_A\,$  of 
\eqref{eq:clustermoments} and of the integrated density states remain
as above, and we maintain the same notation, adding the superscript ``$\ed$''.

For the random walk $(Z_n)$ with law $\mu$ on $\GG\,$, we have
almost surely $[Z_{n-1},Z_n] \in E$ for all $n$.
The following is now a simple exercise.

\begin{lem}\label{lem:bondreturn} The expected $n$-step return
probabilities of the random walk restricted to the bond percolation
cluster are
$$
\Ex_{\psf}^{\ed}\bigl( p_{\omega}^{(n)}(\idG,\idG) \bigr) =
\Ex\bigl(\psf^{|\{ [Z_0, Z_1], \dots, [Z_{n-1},Z_n] \}|}\cdot 
\delta_{\idG}(Z_n)\bigr)\,.
$$
\end{lem}

On the right hand side, expectation refers again to the probability
space underlying the random walk $(Z_n)$.

\subsection{Enlightened edges}
We now consider the situation where a lamp is located on each of the 
\emph{edges} of $X(\GG,S)$ instead of the vertices. When the lamplighter makes
a step from some $x \in \GG$ to a neighbour $y$ (which occurs with probability
$\mu(x^{-1}y)$), then (s)he modifies the
state of the lamp on $\ef = [y,x]$ at random. The state space of this process is
$\CC^{\ed} \rtimes \GG$, where $\CC^{\ed} = \bigoplus_{\Ef} \HH$ consists
of finitely supported configurations 
$\eta: \Ef \to \HH\,$. In this context, the unit element $\uno$ of the group 
$\CC^{\ed}$ is of course given by $\uno(\ef) = \idH$
for all $\ef \in \Ef$. The left action of $\GG$ on $\CC^{\ed}$ is induced by 
the left action on $\Ef$, where $[x,y] \mapsto [gx,gy]$ for $g \in \GG\,$. 
Note that this action is not transitive, but has finitely many orbits.
Of course, $\GG$ is embedded as a subgroup via $g \mapsto (\uno,g)$.
For $\ef \in \Ef$, we define again the measures
$$
\nu_{\ef}(\eta,g) = \begin{cases} 
    1/|\HH|\,,&\text{if}\;g = \idG \;\text{and}\;\supp(\eta) \subset \{\ef\}\,,\\
    0\,,&\text{otherwise;}                              
                    \end{cases}
\qquad\quad \nnu_{\ef} = \delta_{(\uno,\idG)} - \nu_{\ef}\,.
$$		    
Then the (symmetric) law of the new random walk is
$\wt\mu_{\ed} = \sum_{s \in S} \mu(s)\,\cdot\, \nu_{[\idG,s]} * \delta_s\,$,
that is,
\begin{equation}\label{eq:mubond}
\wt\mu_{\ed}(\eta,g) = \begin{cases} 
    \mu(g)/|\HH|\,,&\text{if}\;g \in S \;\text{and}\;\supp(\eta) 
    \subset \{[\idG,g]\}\,,\\
    0\,,&\text{otherwise.}                        
                \end{cases}
\end{equation}

\begin{lem}\label{lem:return-edge}
The return probabilities of the random walk with edge enlightenment 
are
$$
\wt\mu_{\ed}^{(n)}(\uno,\idG) = \Ex\bigl(\psf^{|\{ [Z_0, Z_1], \dots, [Z_{n-1},Z_n] \}|}\cdot 
\delta_{\idG}(Z_n)\bigr)\,.
$$
\end{lem}

We leave the proof once more as an exercise. The basic and well known 
principle (here as well as for the switch-walk-switch random walk) is the 
following.
Up to time $n$, the lamp states can be modified only on those edges 
which the lamplighter crosses, that is, $[Z_0, Z_1], \dots, [Z_{n-1},Z_n]$.
in order to have all lamps switched off at time $n$, the lamplighter must choose
to switch each lamp off at the last visit. This is done with
probability $\psf^{|\{ [Z_0, Z_1], \dots, [Z_{n-1},Z_n] \}|}\,$, after which one
has to average over all possibilities.

Lemmas \ref{lem:bondreturn} and \ref{lem:return-edge} yield the bond-analogue
of Theorem \ref{thm:spec-meas}.

\begin{cor}\label{cor:bond-meas}
The Plancherel measure $\msf^{\ed}$ of the of the edge enlightening random 
walk with law $\wt\mu_{\ed}$ on $ \CC^{\ed} \rtimes\GG$ and the integrated density of states
of bond percolation with parameter $\psf=1/|\HH|$ on $X(\GG,S)$
are related by
$$
\msf^{\ed}(B) = \Ex_p \bigl( \nsf_{\omega}^{\ed}(B) \bigr) \quad 
\text{for every Borel set}\;\;
B \subset \R\,.
$$
\end{cor}

\subsection{Point spectrum}

If we have a finite set $F = \{ [x_1,y_1], \dots, [x_k,y_k] \} \subset \Ef$,
then we define
$$
\nu_F = \nu_{[x_1,y_1]} * \dots \nu_{[x_k,y_k]} \AND
\nnu_F = \nnu_{[x_1,y_1]} * \dots \nnu_{[x_k,y_k]}\,.
$$
The analogue of \eqref{eq:nuCdA} is the signed measure
\begin{equation}\label{eq:nuEdA}
\nu_{E(A),\partial A} = \nu_{E(A)} * \nnu_{\partial A}
\end{equation}
on $ \CC^{\ed}\rtimes \GG\,$, where $A \in \Ccal^{\ed}$. 
We write $Q_A^{\ed}$ for the associated right convolution operator.
It is a projection of $\ell^2( \CC^{\ed}\rtimes \GG)$ onto its image $Q_A^{\ed}\ell^2$.

The analog of Proposition \ref{pro:decompose} is the following.

\begin{pro}\label{pro:bond-decompose} \emph{(a)} The projections $Q_A^{\ed}\,$, 
$A \in \Ccal^{\ed}$, are mutually orthogonal.\\[3pt]
\emph{(b)} We have 
$$
\ell^2( \CC^{\ed}\rtimes \GG) = \bigoplus_{A \in \Ccal^{\ed}} Q_A^{\ed}\ell^2
\quad \text{(closed direct sum), or equivalently,} \quad
\delta_{(\uno, \idG)} = \sum_{A \in \Ccal^{\ed}} \nu_{E(A),\partial A}
$$
if and only if $\theta(\psf)=0$ for bond percolation on $X(\GG,S)$ with 
parameter $\psf=1/|\HH|\,$.
\end{pro}

The proof is basically the same as that of  Proposition \ref{pro:decompose},
observing that for $A \in \Ccal^{\ed}\,$,
$$
\nu_{E(A),\partial A}(\uno,\idG) = \psf^{|E(A)|} (1 - \psf)^{|\partial A|}
= \Prob_{\psf}[C^{\ed}(\idG) = A]\,.
$$
Next, for $x \in A$ and $s \in S$
$$
\nu_{E(A),\partial A} * \delta_x * \nu_{[e,s]} * \delta_s 
= \nu_{E(A),\partial A} * \nu_{[x,xs]} * \delta{xs} =
\begin{cases} \nu_{E(A),\partial A}* \delta_{xs}\,,&\text{if}\; [x,xs] \in A\,,\\
              0\,,&\text{if}\; [x,xs] \in \partial A\,.
\end{cases}
$$	
This yields the analogue of Lemma \ref{lem:invariant}:

\begin{lem}\label{lem:bond-invariant} Let $A \in \Ccal^{\ed}$.  
Then for $f: \GG \to \R$ with  $\supp(f) \subset A$, 
$$
\nu_{E(A),\partial A} * f * \wt \mu_{\ed} = \nu_{E(A),\partial A} * P_Af\,,
$$
and the linear span 
$$
\Bigl\{ \nu_{E(A),\partial A} * f \;:\; f \in \R^{\GG}\,,\; \supp(f) \subset A \Bigr\}
$$
is mapped into itself under right convolution with $\wt \mu_{\ed}$. 
\end{lem}

For each $A \in \Ccal^{\ed}$, we can now choose an orthonormal system 
$f_{A,x}$ of right eigenfunctions of $P_A\,$ with associated
eigenvalues $\la_{A,x}\,$, $x \in A$. Then we define
\begin{equation}\label{eq:SAx-bond}
\sigma_{A,x}^{\ed} = \nu_{E(A),\partial A} * f_{A,x}\,, \AND S_{A,x}^{\ed}(F) 
= F * \sigma_{A,x}\,, 
\end{equation}
$F \in \ell^2( \CC^{\ed}\rtimes \GG)\,$.
Then $\sigma_{A,x}^{\ed} * \wt \mu_{\ed} = 
\la_{A,x} \cdot \sigma_{A,x}^{\ed}\,$, and the images of $S_{A,x}^{\ed}\,$,
$x \in A$, are eigenspaces of $\Rsf_{\wt \mu_{\ed}}$ that span  
$Q_A^{\ed}\ell^2\,$.
The bond-variants of theorems \ref{thm:pointspec} and \ref{thm:pure}
follow.

\begin{cor}\label{cor:bond-pointsp}
\emph{(I)} The point spectrum of the right convolution operator on 
$\ell^2( \CC^{\ed}\rtimes \GG)$ by the measure $\wt\mu_{\ed}$ of 
\eqref{eq:mubond} comprises the set
$$
\Lambda^{\ed} = 
\bigcup \bigl\{ \spec(P_A) : A \in \Ccal^{\ed} \bigr\}\,.
$$
For each eigenvalue $\lambda \in \Lambda^{\ed}$, the eigenspace contains 
infinitely many linearly independent eigenfunctions with finite support.\\[3pt]
\emph{(II)} If $\psf=1/|\HH|$ is such that bond percolation
on $X(\GG,S)$ satisfies $\theta(\psf) = 0$, that is,
$$
\Prob_{\psf}[C^{\ed}(\idG)\;\text{is finite}] = 1\,,
$$
then $\spec(\Rsf_{\wt\mu_{\ed}})$ is pure point; it is the closure of 
$\Lambda^{\ed}$. 
There is a complete orthonormal system in $\ell^2( \CC^{\ed}\rtimes \GG)$ 
consisting
of finitely supported eigenfunctions associated with the eigenvalues in 
$\Lambda^{\ed}$.
\end{cor}

The diagonalization of the convolution operator also follows the 
same lines as above, and we omit the details.

\section{Final remarks}\label{sec:final}

\subsection{Arbitrary $\psf$} We can start with $\GG$ and a symmetric probability measure
$\mu$ as above, and consider site (or bond) percolation with an arbitrary
parameter $\psf \in (0\,,\,1)$. Then we have the spectrum and associated 
integrated density of states of the (substochastic) random transition 
operator on the cluster $C_{\omega}(\id)$. Can it always be described as
the Plancherel measure of a deterministic convolution operator $\Rsf_{\wt \mu}$ 
on a wreath product $\HH \wr \GG$~? If $\psf = 1/N$ then we can take 
$\wt \mu = \nu * \mu * \nu$, where $\nu$ is equidistribution on
a finite group $\HH$ of order $N$.

For other values of $\psf$, we can proceed as in \cite{DiSc}. We need a group 
$\HH$ and a (signed) measure $\nu$ on $\HH$ which satisfies
$$
\nu(\idH) = \psf \AND  \nu * \nu = \nu\,.
$$
This means that the right convolution operator $\Rsf_\nu$ is a projection with von
Neumann trace $\psf$.
Once we have such a measure $\nu$, we can again define the (signed) measure
$\wt \mu = \nu * \mu * \nu$, where $\mu$ and $\nu$ are now considered as 
measures on the wreath product $\HH \wr \GG$ via the
natural embeddings of $\GG$ and $\HH$. Apart from
the fact that $\wt \mu$ is not a probability measure, all the above results and 
computations remain unchanged. However the eigenfunctions found in this way are 
no more finitely supported, unless the range of the projection $\Rsf_\nu$
has a basis of finitely supported functions.

For example, we can always take $\HH = \Z$ and the measure $\nu$ whose
Fourier transform (characteristic function) is 
$\mathbf{1}_{[-\pi \psf\,,\,\pi \psf]}\,$. That is,
$$
\nu(0) = \psf\,, \AND \nu(k) = \frac{\sin(k\, \pi\, \psf)}{k \,\pi}\,,\;
(k \ne 0)\,.
$$ 

\subsection{Continuous spectrum} We do not know what happens with the rest 
of the spectrum in
the case when $\theta(\psf) >0$ for site percolation on $X(\GG,S)$. 
As Theorem \ref{thm:pointspec} and Corollary \ref{cor:spec} show,
we always have a big point spectrum, but by Proposition \ref{pro:decompose}, 
a part of the Plancherel measure of $\Rsf_{\wt\mu}$ is still missing.
The situation for bond percolation, resp. enlightened edges, is the same.

As a matter of fact, it is an open question whether the integrated density
of states of percolation on Cayley graphs has a continuous part.
It is known that the rest of the spectrum (corresponding to infinite clusters)
is not purely continuous, see e.g.\ \textsc{Kirkpatrick and Eggarter}~\cite{KiEg}
and \textsc{Veseli\'c}~\cite{Ve}, but it is expected to have
a continuous part at least in certain cases.
[We acknowledge
feedback of Peter M\"uller (G\"ottingen) on this question, which reached us
via Florian Sobieczky.]
Maybe the translation of this spectral problem from a random operator
to a deterministic one on the lamplighter group can help to find an answer.

\subsection{Generalizations} We have resisted the temptation to build up a generalized
setting in which both models (switch-walk switch and enlightened edges,
or site and bond percolation, respectively), as well as further variants,
arise as special cases. Our aim was to clarify the relation between the spectra
of lamplighter walks and site \& bond percolation in the most basic context.

Bond percolation on $X(\GG,S)$ is of course equivalent with site percolation
on the line graph of $X(\GG,S)$. The vertex set of the latter is the
edge set of the original graph, and two such edges are now neighbours 
if they have a common end point. Similarly, the random walk with
enlightened edges is up to a few adaptations equivalent with a 
switch-walk-switch lamplighter walk over the line graph. Further types of 
lamplighter random walks can also be related with percolation
(e.g. oriented percolation).

Most of what is exhibited here can be extended to vertex-transitive graphs,
and an adaptation of a large part to general locally finite graphs will 
also not be too hard.

We also remark that the convolution operator on $\HH \wr \Z$ considered
by {\sc Grigorchuk and \.Zuk}~\cite{GrZu} and 
{\sc Dicks and Schick}~\cite{DiSc} is not exactly the same as the 
switch-walk-switch lamplighter walk. Informally, it is the model 
``switch first and then walk to the right, or walk first to the left and 
then switch at the arrival
point''. By shifting the lamps to the edges, one sees that this model
is completely equivalent to the one with enlightened edges on $\Z$. 

\subsection{Finite Groups} While we usually have in mind that $\GG$ is an 
\emph{infinite,} finitely generated group, this is not really relevant. 
If $\GG$ is finite, then the method of 
\S\ref{sec:proofs}.\ref{ssec:diagonalizationtorsion} applies for diagonalizing 
the transition matrix of a symmetric switch-walk-switch random walk on 
$\HH \wr \GG$, where both groups are finite. Explicit computations can 
be performed for example when $X(\GG,S)$ is a cycle or a complete graph.
A quite different approach to deal with this class of examples
is used in recent work of {\sc Scarabotti and Tolli}~\cite{ScTo}.

\end{document}